\input amstex
\magnification=\magstep0
\documentstyle{amsppt}
\pagewidth{6.2in}
\topmatter
\title
Play time with determinants
\endtitle
\author
Tewodros Amdeberhan and Shalosh B. Ekhad
\endauthor
\endtopmatter
\def\({\left(}
\def\){\right)}

\document
\noindent
Computing determinants is invariably in demand from all sorts of mathematical domains and needs neither cynical advocacy nor does it lack motivation. In this semi-expository article we wish to illustrate a few techniques that should receive ample attention due to their simplicity, elegance and effective nature.
\smallskip
\noindent
\bf 0 Glossary of notations and conventions \rm
\smallskip
\noindent
Throughout, almost all determinants have indices $(i,j)$ ranging in $\{0,1,\dots,n-1\}$. Where there is no confusion, we omit to write them.
\smallskip
\noindent
$m!!=1!2!\cdots m!$ and $(x)_i=x(x+1)\cdots(x+i-1)$ is the Pochhammer symbol.
\smallskip
\noindent
$\binom{a}b=0$ if either $b<0$ or $b>a$.
\smallskip
\noindent
$\widehat{Ct_i}$ is the constant-term operator that extracts the coefficient of $x_i^0$. Also $\widehat{Ct}=\prod_{i=0}^{n-1}\widehat{Ct_i}$.
\smallskip
\noindent
$S_n$ is the symmetric group of permutations on $n$ letters.
\smallskip
\noindent
$\bold{x}_n=(x_0,x_1,\dots,x_{n-1})$ and $\bold{X}_n=x_0x_1\cdots x_{n-1}$.
\smallskip
\noindent
$V(\bold{x}_n)$ is the Vandermionde determinant $\prod_{j>i}(x_j-x_i)$.
\smallskip
\noindent
$[n]_q=\frac{1-q^n}{1-q}$ and $[n]!_q=[1]_q[2]_q\cdots[n]_q$ 
\smallskip
\noindent
The Gaussian polynomials $\binom{n}k_q=\frac{[n]!_q}{[k]!_q[n-k]!_q}$, for $0\leq k\leq n$; $\binom{n}k_q=0$ for $k<0$ or $k>n$.

\bigskip
\noindent
\bf 1 The Method of Condensation \rm
\smallskip
\noindent
The idea originated with Charles L. Dodgson and the iterative form is propelled by Zeilberger [7]. Here are some instances of its application.
\smallskip
\noindent
\bf 1.1 Example. \rm Consider the $n\times n$ matrix of entries $\binom{2i}j$. Let us compute the determinant. The first step is to generalize it (introduce two additional parameters) as $Z_n(a,b)=\det\left[\binom{2i+2a}{j+b}\right]$. The next step is an automated guess for its determinantal value
$$2^{\binom{n}2}\prod_{i=0}^{n-1}\frac{(2i+2a)!i!}{(i+b)!(2i+2a-b)!}=2^{\binom{n}2}\prod_{i=0}^{n-1}\binom{2i+2a}b\binom{i+b}b^{-1}.\tag1.1.1$$
\bf Proof. \rm According to [7], the following is always true for any determinant (when $i\rightarrow i+a; j\rightarrow j+b$)
$$Z_n(a,b)=\frac{Z_{n-1}(a,b)Z_{n-1}(a+1,b+1)-Z_{n-1}(a+1,b)Z_{n-1}(a,b+1)}{Z_{n-2}(a+1,b+1)}.$$
It remains to verify that (1.1.1) also satisfies this non-linear recurrence. Symbolic softwares are naturally adept at this task. The proof will then be complete once two initial cases are checked, say for $n=1$ and $n=2$. $\square$
\smallskip
\noindent
\bf 1.2 Example. \rm In particular, setting $a=b=0$ implies that 
$\det\left[\binom{2i}j\right]=2^{\binom{n}2}$. The Online Encyclopedia of Integer Sequences (OEIS) associates $2^{\binom{n}2}$ with the number of $n\times n$ binary matrices with no row sum greater than $n-1$, or the number of labeled $n$-colorable graphs on $n$ nodes.
\smallskip
\noindent
We introduce a generalization of formula (1.1.1) that is provable by Condensation.
$$\det\left[\binom{ri+x}{j+y}\right]=
r^{\binom{n}2}\prod_{i=0}^{n-1}\binom{ri+x}{y}\binom{i+y}y^{-1}
.\tag1.2.1$$
Once again, the specialization $x=y=0$ yields a pretty identity:
$\det\left[\binom{ri}j\right]=r^{\binom{n}2}$.
Equation (1.2.1) also shows something interesting: the following determinant is independent of $x$.
$$\det\left[\binom{qi+x}j\right]=q^{\binom{n}2}.$$
\bf 1.3 Example. \rm The number of plane partitions that fit in an $a\times b\times c$ box is enumerated by
$$D_n(a,b):=\det\left[\binom{i+j+a+b}{i+a}\right]_{i,j=0}^{c-1}
=\prod_{k=0}^{c-1}\prod_{j=0}^{b-1}\prod_{i=0}^{a-1}\frac{i+j+k+2}{i+j+k+1}.$$
Again, this formula is provable by the Condensation method. This determinantal representation is recorded here due to the apparent symmetry of the matrix on the left-hand side. 
\bigskip
\noindent
\bf 2 Buy one, get three free! \rm
\smallskip
\noindent
In the last example, set $a=b=0$, and write $D_n$ instead for $D_n(0,0)$.
We will calculate $D_n$ in three different ways, each leading to a different interpretation.
\bigskip
\noindent
\bf 2.1 \`A la Andrews \rm
\smallskip
\noindent
The present method has been utilized by George Andrews [2] in his proof of the Mills-Robbins-Rumsey determinant [6]
$$A_n:=\det\left[\binom{\mu+j+i}{2j-i}\right]=2^{-n}\prod_{i=0}^{n-1}\frac{(\mu+2i+2)_i(\frac12\mu+2i+\frac32)_{i-1}}{(i)_i(\frac12\mu+i+\frac32)_{i-1}}.\tag2.1.1$$
The main essence was to exhibit a triangular matrix $T_n$ so that $A_nT_n$ is triangular too. In the same spirit, let us implement the Vandermonde-Chu identity 
$$\binom{i+j}i=\sum_k\binom{i}k\binom{j}{i-k},$$ 
which follows from an elementary combinatorial argument. Therefore $$D_n=\det\binom{i+j}i=\det\binom{i}k\cdot\det\binom{k}j=1,$$ 
since in the latter the matrices are triangular with $1^{\prime}$s on the main diagonal.

\pagebreak
\bigskip
\noindent
\bf 2.2 Constant-term identities \rm
\smallskip
\noindent
This method goes back, to at least, Zeilberger$^{\prime}$s paper [8]. For the current application, begin as follows:
$$\binom{i+j}i=\widehat{Ct_i}\left(x_i^i(1+x_i^{-1})^{i+j}\right)=\widehat{Ct_i}(1+x_i)^i(1+x_i^{-1})^j.$$ 
Therefore
$$\align
D_n
&=\widehat{Ct}\prod_{i=0}^{n-1}(1+x_i)^i\det\left[(1+x_i^{-1})^j\right]\\
&=\widehat{Ct}\prod_{i=0}^{n-1}(1+x_i)^i\prod_{j>i}\left(x_j^{-1}-x_i^{-1}\right)\\
&=\frac1{n!}\widehat{Ct}\prod_{j>i}\left(x_j-x_i\right)\prod_{j>i}\left(x_j^{-1}-x_i^{-1}\right)\\
&=\frac1{n!}\widehat{Ct}\frac{(-1)^{\binom{n}2}}{\prod_{i=0}^{n-1}x_i^{n-1}}\prod_{j>i}\left(x_j-x_i\right)^2.
\endalign$$
Since $\prod_{j>i}(x_j-x_i)$ is the Vandermonde determinant $V(\bold{x})$, we gather the coefficient of $\bold{X}_n^{n-1}$ in $V^2(\bold{x}_n)$ to be $(-1)^{\binom{n}2}n!$. This, of course, reproves a trivial special case of Dyson$^{\prime}$s conjecture that the constant term in $\prod_{j\neq i}\left(1-\frac{x_i}{x_j}\right)$ equals $n!$.
\smallskip
\noindent
\bf 2.2.1 Remark. \rm The \it even \rm powers (in particular, the square) of the Vandermonde determinant play a crucial role in the quantum Hall effect phenomena via Laughlin$^{\prime}$s wave function ansatz. Determining the coefficients of their expansion, as symmetric functions, in terms of Schur functions has triggered considerable interest among physicists (see [3] and references therein). At the moment of this writing, extracting these numbers in general is an open problem.
\bigskip
\noindent
\bf 2.3 Multilinearity of the determinant \rm
\smallskip
\noindent
Define the determinant $B_n=\det[(i+j)!]$. For the third computation, watch an online video of a lively lecture delivered by Zeilberger [10]. For the purpose at hand, let us make use of Euler$^{\prime}$s formula $m!=\int_0^{\infty}t^me^{-t}dt$. Thus
$$\align
B_n&=\det\left[\int_0^{\infty}x_j^{i+j}e^{-x_j}dx_j\right]\\
&=\int_{\Bbb{R}_{+}^n}\left(e^{-\sum x_j} V(\bold{x}_n)
\prod_{i=0}^{n-1}x_j^j\right)d\bold{X}_n\\
&=\frac1{n!}\sum_{\pi\in S_n}\int_{\Bbb{R}_{+}^n}\left(e^{-\sum x_j} V(\bold{x}_n)
(-1)^{\pi}\prod_{i=0}^{n-1}x_{j}^{\pi(j)}\right)d\bold{X}_n\\
&=\frac1{n!}\int_{\Bbb{R}_{+}^n}e^{-\sum x_j} V^2(\bold{x}_n)
d\bold{X}_n.\endalign$$
\smallskip
\noindent
Observe that $B_n=\det[i!j!\binom{i+j}i]=(n-1)!!^2D_n=(n-1)!!^2$. Hence, we obtain the Selberg-type integral
$$\int_{\Bbb{R}_{+}^n}e^{-\sum x_j} V^2(\bold{x}_n)
d\bold{X}_n=n!!(n-1)!!.$$
\bf 2.4 Remark. \rm In summary, the above ideas combined with Dodgson$^{\prime}$s Condensation suggest simple techniques that are efficient in proving more general evaluations, such as
$$\det\left[\binom{i+j+\alpha}{i+\beta}\right],\qquad
\widehat{Ct}\prod_{j\neq i}\left(1-\frac{x_i}{x_j}\right)^{\alpha},\qquad
\int_{\Bbb{R}_{+}^n}\bold{X}_n^{\alpha}e^{-\sum x_j} V^{2\beta}(\bold{x}_n)d\bold{X}_n.$$
\bigskip
\noindent
\bf 3 More Examples \rm
\smallskip
\noindent  
The Delannoy numbers $D(i,j)$ count walks from $(0,0)$ to $(i,j)$ composed of unit steps East, North, and North-East. These are given by $\sum_{k\geq0}\binom{i}k\binom{j}k2^k$. 
\bigskip
\noindent
\bf 3.1 Example. \rm The below determinant was conjectured by Bacher and Krattenthaler [5] provided a proof. We offer a simpler argument using Constant-term identities. For alternative automatic proof see Koutschan [4].
$$K_n:=\det\left[\sum_{k=0}^{\min(2i,2j)}\binom{2i}k\binom{2j}k2^k\right]=4^{\binom{n}2}\prod_{i=0}^{n-1}\frac{i!^2(4i)!}{(2i)!^3}.$$
\bf Proof. \rm Start with $\sum_{k\geq0}\binom{2i}k\binom{2j}k2^k=\widehat{Ct_i}(1+2x_i)^{2i}(1+x_i^{-1})^{2j}$. Then
$$\align 
K_n&=\widehat{Ct}\prod_{i=0}^{n-1}(1+2x_i)^{2i}\det\left[(1+x_i^{-1})^{2j}\right]\\
&=\frac{4^{\binom{n}2}}{n!}\widehat{Ct}\prod_{j>i}(x_j+x_j^2-x_i-x_i^2)\det\left[(1+x_i^{-1})^{2j}\right]\\
&=4^{\binom{n}2}\widehat{Ct}\prod_{i=0}^{n-1}(x_i+x_i^2)^i\det\left[(1+x_i^{-1})^{2j}\right]\\
&=4^{\binom{n}2}\det\left[\widehat{Ct_i}(x_i+x_i^2)^i(1+x_i^{-1})^{2j}\right]\\
&=4^{\binom{n}2}\det\left[\widehat{Ct_i}\left(x^{i-2j}(1+x_i)^{2j+i}\right)\right]\\
&=4^{\binom{n}2}\det\left[\binom{2j+i}{2j-i}\right].\endalign$$
Formula (1) in [1], with $x=y, a=2$, computes the last determinant to produce the desired value. $\square$
\smallskip
\noindent
\bf 3.2 Example. \rm Some play time with variants of Example 3.1 lead us to conjecture the determinant $L_n$ given below. The authors are grateful to Tiago Dinnis Da Fonseca (Univ. Montreal) who allowed his proof to be included in this article. It is rather amusing how $L_n$ turns out to be equivalent to (2.1.1). By association, its evaluation is as simple or as complicated. For alternative automatic proof see Koutschan [4]. $$L_n:=\det\left[\sum_{k=0}^{\min(2i,2j)}\binom{2i}k\binom{2j}k4^k\right]
=16^{\binom{n}2}\prod_{i=0}^{n-1}\frac{(2i)!(6i)!(3i+1)}{(4i)!^2(4i+1)}.$$
\bf Proof. \rm It is evident that $\sum_{k\geq0}\binom{2i}k\binom{2j}k4^k=\widehat{Ct_i}(1+2x_i)^{2i}(1+2x_i^{-1})^{2j}$. Then
$$\align 
L_n&=\widehat{Ct}\prod_{i=0}^{n-1}(1+2x_i)^{2i}\det\left[(1+2x_i^{-1})^{2j}\right]\\
&=\frac{4^{\binom{n}2}}{n!}\widehat{Ct}\prod_{j>i}(x_j+x_j^2-x_i-x_i^2)\det\left[(1+2x_i^{-1})^{2j}\right]\\
&=4^{\binom{n}2}\widehat{Ct}\prod_{i=0}^{n-1}(x_i+x_i^2)^i\det\left[(1+x_i^{-1})^{2j}\right]\\
&=4^{\binom{n}2}\widehat{Ct}\prod_{i=0}^{n-1}(x_i+x_i^2)^i4^{\binom{n}2}\prod_{j>i}(x_j^{-1}+x_j^{-2}-x_i^{-1}-x_i^{-2})\\
&=16^{\binom{n}2}\widehat{Ct}\prod_{i=0}^{n-1}(x_i+x_i^2)^i\det\left[(x_i^{-1}+x_i^{-2})^j\right]\\
&=16^{\binom{n}2}\det\left[\widehat{Ct_i}(x_i+x_i^2)^i(x_i^{-1}+x_i^{-2})^j\right]\\
&=16^{\binom{n}2}\det\left[Ct_i\left(x^{i-2j}(1+x_i)^{i+j}\right)\right]\\
&=16^{\binom{n}2}\det\left[\binom{j+i}{2j-i}\right].\endalign$$
Alas! this is the familiar Mills-Robbins-Rumsey determinant (2.1.1), specialized at $\mu=0$,
$$\det\left[\binom{j+i}{2j-i}\right]=2^{-n}\prod_{i=0}^{n-1}\frac{(2i+2)_i(2i+\frac32)_{i-1}}{(i)_i(i+\frac32)_{i-1}}.$$
Using $(y)_i=\frac{(y+i-1)}{(y-1)!}$, the right-hand side routinely transforms to the required product expression. $\square$
\bigskip
\noindent
\bf 4 The Holonomic Ansatz \rm
\smallskip
\noindent
Once more, this automatic method for the evaluation of determinants is a creation of Zeilberger [9]. We are grateful to Christoph Koutschan (RISC-Linz) for allowing his proofs, which exploit the present method of holonomic ansatz, to be made available [4]. So far, this is the only technique providing a proof to the following ex-conjecture.
\bigskip
\noindent 
\bf 4.1 Example. \rm We have the determinantal evaluation
$$\align
\det\left[\sum_{k=0}^{\min(3i,3j)}\binom{3i}k\binom{3j}k3^k\right]
&=2^{8\binom{n}2}\prod_{i=0}^{n-1}
\frac{(\frac7{12})_i(\frac1{12})_i(\frac54)_i(\frac34)_i}{(\frac76)_i(\frac16)_i(\frac23)_i(\frac23)_i}\\
&=6^{2\binom{n}2}\prod_{i=1}^{n-1}\prod_{j=1}^i
\frac{(12j-5)(12j-11)(4j+1)(4j-1)}{(6j+1)(6j-5)(3j-1)(3j-1)}.\endalign$$

\pagebreak
\bigskip
\noindent
\bf 5 Quantum Analogues \rm
\smallskip
\noindent
After the initial release of the present note and its conjectures, the authors received instant feedbacks. They are grateful to Johann Cigler (Univ. Vienna) for the permission to incorporate his conjectures and comments on certain $q$-analogues. Although the proofs might be carried out by employing the methods from the earlier sections, we are not pursuing the details. The reader is invited to do so. 
\bigskip
\noindent
\bf 5.1 Example. \rm A $q$-analogue of formula (1.2.1) is given by
$$\det\left[\binom{ri+x}{j+y}_qq^{i(y-x)}\right]
=[r]_q^{\binom{n}2}\prod_{i=0}^{n-1}q^{(r-1)\binom{i}2}\frac{[i]!_{q^r}\binom{ri+x}y_q}{[i]!_q\binom{i+y}y_q}.$$
\bf 5.2 Example. \rm The determinants
$$\det\left[\sum_{k=0}^{\min(i,j)}\binom{ri}{k}_q\binom{sj}k_qz^k\right]
=z^{\binom{n}2}\prod_{i=0}^{n-1}\binom{ri}i_q\binom{si}i_q, \qquad\text{and}$$
$$\det\left[\sum_{k=0}^{\min(i,j)}\binom{i}{k}_q\binom{j}k_q\prod_{l=1}^k(1+q^l)\right]
=\prod_{i=0}^{n-1}(1+q^i)^{n-i}$$
both follow from the general formula
$$\det\left[\sum_{k=0}^{\min(i,j)}a(i,k)b(j,k)f(k)\right]
=\prod_{i=0}^{n-1}a(i,i)b(i,i)f(i).\tag5.2.1$$
The identity (5.2.1) results from an immediate application of the method outlined in Section 2.1. 
\bigskip
\noindent
\bf 5.3 Example. \rm 
The determinant
$$\det\left[\sum_{k=0}^{\min(ri,j)}\binom{ri}k_q\binom{rj}k_q\right]
=q^{(r-1)\binom{n}3}\prod_{j=1}^{n}[r]_{q^j}^{n-j}\prod_{i=0}^{n-1}\binom{ri}i_q$$
follows from
$$\det\left[\sum_{k=0}^{\min(ri,j)}\binom{ri}k_qf(j,k)\right]
=\det\left[\sum_{k=0}^{\min(ri,j)}\binom{ri}k_q\right]_{i,j=0}^{n-1}
\prod_{i=0}^{n-1}f(i,i).$$
Therefore, it suffices to prove that
$$\det\left[\sum_{k=0}^{\min(ri,j)}\binom{ri}k_q\right]=
q^{(r-1)\binom{n}3}\prod_{j=1}^{n}[r]_{q^j}^{n-j}.$$

\enddocument